\documentclass{article}
\usepackage{amsfonts, amssymb, latexsym, amscd}
\newtheorem{thm}{Theorem}[section]
\newtheorem{cor}[thm]{Corollary}
\newtheorem{lem}[thm]{Lemma}
\newtheorem{prp}[thm]{Proposition}

\def\qed{{\hfill$\Box$}}
\def\T{{\rm \Delta}}
\def\g{{\cal G}}
\def\gg{{\cal G'}}
\def\ggg{{\cal G^*}}

\begin{document}

\title{\bf Spectral characterization of
\\
mixed extensions of small graphs}

\author{Willem H. Haemers\thanks{e-mail haemers@uvt.nl}
\\
{\it\small Department of Econometrics and Operations Research,}
\\
{\it\small Tilburg University, The Netherlands}
}
\date{}

\maketitle

\begin{abstract}
\noindent
A mixed extension of a graph $G$ is a graph $H$ obtained from $G$ by replacing each vertex of
$G$ by a clique or a coclique, where vertices of $H$ coming from different
vertices of $G$ are adjacent if and only if the original vertices are adjacent
in $G$.
If $G$ has no more than three vertices, $H$ has all but at most three adjacency
eigenvalues equal to $0$ or $-1$.
In this paper we consider the converse problem, and determine the class $\g$
of all graphs with at most three eigenvalues unequal to $0$ and $-1$.
Ignoring isolated vertices, we find that $\g$ consists of all mixed extensions of
graphs on at most three vertices together with some particular mixed extensions of the
paths $P_4$ and $P_5$.
\\[5pt]
{Keywords:}~graph spectrum, spectral characterization.
AMS subject classification:~05C50.
\end{abstract}

\section{Introduction}

We determine the class $\g$ of graphs (simple and undirected) with
the property that all but at most three eigenvalues of the adjacency matrix are equal to $-1$ or $0$.
The research is motivated by the question for which values of $p$ and $q$ the pineapple graph
$K_p^q$ (we will give the definition below) is determined by the spectrum of the adjacency matrix.
The pineapple graph has the above mentioned property, and a partial answer to the question is given in \cite{hst}.
The classification of $\g$ can be an important step towards a complete answer,
and moreover, it may lead to the determination of all graphs in $\g$ determined by the spectrum.
A comparable approach has been successful for graphs with all but two eigenvalues equal to $\pm 1$ (see~\cite{chvw}),
and $-2$ or $0$ (see~\cite{chv}).

Deleting or adding an isolated vertex does not effect the mentioned property, therefore we can
restrict to the class $\gg$ of graphs in $\g$ with no isolated vertex.
It will turn out that $\gg$ can be described in terms of mixed extensions, which are introduced
in the next section.
\\

We assume familiarity with basic results from linear algebra and graph spectra.
Our main tool is the interlacing theorem presented below (we refer to \cite{cds}, and \cite{bh}
for this and other results on graph spectra).
\begin{thm}
Let $H$ be a graph of order $m$ with (adjacency) eigenvalues $\lambda_1\geq\ldots\geq\lambda_m$,
and let $G$ be an induced subgraph of $H$ of order $n$ with eigenvalues $\mu_1\geq\ldots\geq\mu_n$.
Then
\[
\lambda_i\geq\mu_i\geq\lambda_{m-n+i}\ \ \mbox{for}\ i=1,\ldots,n.
\]
\end{thm}
We illustrate the use of this theorem by proving two well-known characterizations,
which are relevant for our approach (the second one is due to Smith~\cite{s}).
\begin{prp}\label{K}
(i)
If a graph $H$ has smallest eigenvalue at least $-1$, then $H$ is the disjoint union of complete graphs.
\\
(ii)
If a graph $H$ has just one positive eigenvalue, then $H$ is a complete multipartite graph,
possibly extended with some isolated vertices.
\end{prp}
{\bf Proof.}
(i)
The path $P_3$ has smallest eigenvalue $-\sqrt{2}<-1$, and therefore is not an induced subgraph of $H$.
Hence each connected component of $H$ is a clique.
\\
(ii) Assume $H$ has no isolated vertices and $G=K_2+K_1$ (the disjoint union of an edge and an isolated vertex) as an induced subgraph.
Then the isolated vertex in $G$ cannot be isolated in $H$, 
therefore $K_2+K_2$, $P_4$, or $K_3$ with a pendant edge is an induced subgraph of $H$, 
however each of these graphs has two positive eigenvalues,
which violates the interlacing inequalities.
Hence $G$ is not an induced subgraph of $H$, which implies that $H$ is a complete multipartite graph.
\qed
\section{Mixed extensions}
Consider a graph $G$ with vertex set $\{1,\ldots,n\}$.
Let $V_1,\ldots,V_n$ be mutually disjoint nonempty finite sets.
We define a graph $H$ with vertex set the union of $V_1,\ldots,V_n$ as follows.
For each $i$, the vertices of $V_i$ are either all mutually adjacent ($V_i$ is a clique),
or all mutually nonadjacent ($V_i$ is a coclique).
When $i\neq j$, a vertex of $V_i$ is adjacent to a vertex of $V_j$ if and only if $i$ and $j$ are adjacent in $G$.
We call $H$ a {\em mixed extension} of $G$.
We represent a mixed extension by an $n$-tuple $(t_1,\ldots,t_n)$ of nonzero integers,
where $t_i>0$ indicates that $V_i$ is a clique of order $t_i$, and $t_i<0$ means that $V_i$ is a coclique of order $-t_i$.
An empty position means that $t_i=\pm 1$. 

For example, for positive $p$ and $q$, the mixed extension of $K_2$ of type $(-p,-q)$ is the complete bipartite graph $K_{p,q}$,
and the mixed extension of $K_2$ of type $(-p,q)$, is the complete multipartite graph $K_{p,1,\ldots,1}$,
also known as the {\em complete split graph} $CS_{p,q}$.
The mixed extension of the path $P_3$ of type $(p,\  ,-q)$ is the {\em pineapple graph} $K_{p+1}^q$.
If $t_i=p>1$ for $i=1,\ldots,n$ we speak of a $p$-{\em clique extension}, and if $t_i=-q<-1$ for $i=1,\ldots,n$
then it is called a {\em $q$-coclique} extension.
If the mixed extension involves only cliques ($t_i\geq -1$), or cocliques ($t_i\leq 1$),
we call it a {\em mixed clique}, or {\em mixed coclique extension}, respectively.

We denote the all-ones matrix by $J$, and the all-zeros matrix by $O$.
If $A$ and $B$ are the adjacency matrices of $G$ and $H$ respectively, then $B$ can be obtained from $A$
by replacing the $i$-th diagonal entry of $A$ by an $|t_i|\times|t_i|$ block $O$ if $t_i<-1$,
and $J-I$ if $t_i>1$.
In addition, an off diagonal $0$ or $1$ of $A$ is replaced by $O$ or $J$, respectively.
The {\em quotient matrix} of $B$ is the $n\times n$ matrix $Q$ where each entry equals the row sum of the corresponding block in $B$.

\begin{lem}
The eigenvalues of $B$ are the eigenvalues of $Q$ together with $m-n$ eigenvalues equal to $0$ or $-1$.
\end{lem}
{\bf Proof.}
The partition is equitable, therefore the spectrum of $Q$ is a sub-multiset of the spectrum of $B$.
The remaining eigenvalues of $B$ do not change when an all-one block $J$ is subtracted from a block of
$B$ (see~\cite{bh}, Section~2.3).
Doing so for all non-zero blocks yields a diagonal matrix with entries $-1$ or $0$.
\qed
\\[5pt]
In particular, if $H$ is a mixed extension of a graph with no more than three vertices, then $H\in\g$.

\section{Easy cases}\label{easy}

\begin{prp}\label{2}
Suppose $H\in\gg$ has all but at most two eigenvalues equal to $-1$ or $0$, then $H$ is one of the following.
\\
(i) $K_m$\ with $m\geq 2$,
\\
(ii) $K_{p}+K_{q}$\ with $p,q\geq 2$,
\\
(iii) $K_{p,q}$\ with $p,q\geq 2$,
\\
(iv) $CS_{p,q}$\ with $p\geq 2$, $q\geq 1$.
\end{prp}
{\bf Proof.}
Clearly $H$ has at least one positive eigenvalue.
If all other eigenvalues are equal to $-1$ or $0$, then $H$ has least eigenvalue $-1$,
hence $H$ is the complete graph $K_m$.
Next assume $H$ has exactly two eigenvalues different from $-1$ and $0$.
If both are positive, then again $-1$ is the smallest eigenvalue
and therefore $H$ is the disjoint union of two cliques with at least two vertices.
If $H$ has one positive eigenvalue and one negative eigenvalue $\neq -1$, then the second largest eigenvalue of $H$
is not positive, so by Proposition~\ref{K}, $H$ is a complete multipartite graph.
The complete $3$-bipartite graph $K_{2,2,1}$ has three eigenvalues not in $[-1,0]$,
therefore, by interlacing, $K_{2,2,1}$ is not an induced subgraph of $H$.
So $H$ is either a complete bipartite graph, or $H=K_{p,1,\ldots,1}=CS_{p,q}$ ($p,q\geq 2$).
\qed
\begin{prp}\label{disconnected}
Suppose $H\in\gg$ is disconnected with exactly three eigenvalues different from $-1$ and $0$, then $H$ is one of the following.
\\
(i) $K_p+K_q+K_r$,\ with $p,q,r\geq 2$,
\\
(ii) $K_p+K_{q,r}$,\ with $p,q,r\geq 2$,
\\
(iii) $K_p+CS_{q,r}$,\ with $p,q\geq 2,\ r\geq 1$.
\end{prp}
{\bf Proof.}
Since $H$ has no isolated vertices, each connected component has all but at most two eigenvalues equal to $-1$ or $0$.
So the possible components are given by Proposition~\ref{2}.
\qed
\begin{prp}\label{3}
Suppose $H$ is a connected graph with exactly three eigenvalues different from $-1$ and $0$, then one of the following holds.
\\
(i) $H=K_{p,q,r}$ with $p,q,r\geq 2$,
\\
(ii) $H$ is a mixed extension of $K_3$ of type $(-p,-q,r)$ with $p,q\geq 2$, $r\geq 1$,
\\
(iii) $H$ has exactly two positive eigenvalues and exactly one eigenvalue less than $-1$.
\end{prp}
{\bf Proof.}
Clearly $H$ has at least one positive eigenvalue.
Assume $H$ has just one positive eigenvalue.
Then $H$ is a complete multipartite graph.
The complete multipartite graph $K_{2,2,2,1}$ has three eigenvalues less than $-1$,
and therefore cannot be an induced subgraph of $H$.
Since $H$ has two eigenvalues less than $-1$, it follows that either
$H=K_{p,q,r}$ with $p,q,r\geq 2$, which is case (i), or $H=K_{p,q,1,\ldots,1}$ with $p,q\geq 2$,
which is case (iii).

If $H$ has three positive eigenvalues, then the smallest eigenvalue equals $-1$, so $H$ is the disjoint union of three complete graphs, which contradicts the assumption.
The remaining possibility is case (iii).
\qed\\

What remains to be determined is the class of connected graphs in $\gg$ with two positive eigenvalues and one
eigenvalue less than $-1$.
We denote this class by $\ggg$.
Fist we make a list of forbidden induced subgraphs.

\section{Forbidden induced subgraphs}\label{forbidden}

\begin{lem}\label{forb}
The following graphs are not induced subgraphs of any graph in $\ggg$:
\\[5pt]
The pentagon $C_5$, the path $P_6$ (and therefore $C_n$ for $n\geq 7$), and
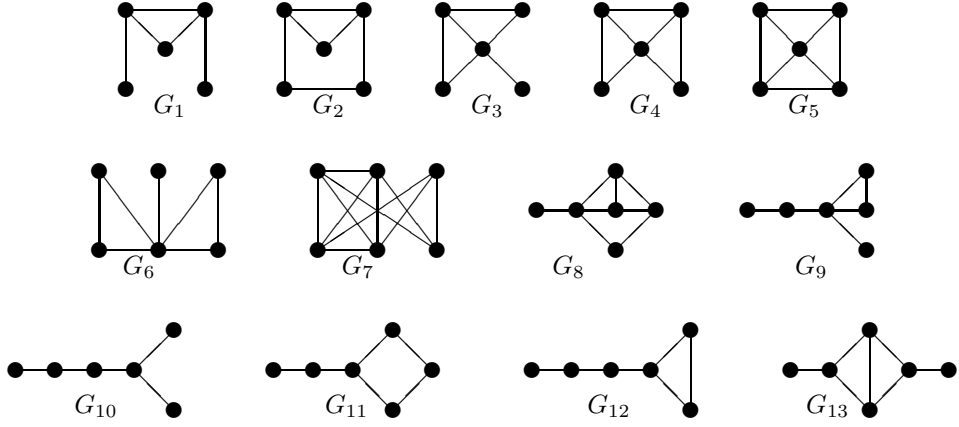
\setlength{\unitlength}{3pt}
\[
\ \ \ \ \ \
\begin{picture}(20,12)(0,0)
\put(5,5){\circle*{2}}
\put(0,0){\circle*{2}}
\put(10,0){\circle*{2}}
\put(10,10){\circle*{2}}
\put(0,10){\circle*{2}}
\put(0,10){\line(1,-1){5}}
\put(0,10){\line(1,0){10}}
\put(10,10){\line(-1,-1){5}}
\put(10,0){\line(0,1){10}}
\put(0,0){\line(0,1){10}}
\put(3.5,-3){$G_1$}
\end{picture}
\begin{picture}(20,12)(0,0)
\put(5,5){\circle*{2}}
\put(0,0){\circle*{2}}
\put(10,0){\circle*{2}}
\put(10,10){\circle*{2}}
\put(0,10){\circle*{2}}
\put(0,10){\line(1,-1){5}}
\put(0,10){\line(1,0){10}}
\put(10,10){\line(-1,-1){5}}
\put(10,0){\line(0,1){10}}
\put(0,0){\line(1,0){10}}
\put(0,0){\line(0,1){10}}
\put(3.5,-3){$G_2$}
\end{picture}
\begin{picture}(20,12)(0,0)
\put(5,5){\circle*{2}}
\put(0,0){\circle*{2}}
\put(10,0){\circle*{2}}
\put(10,10){\circle*{2}}
\put(0,10){\circle*{2}}
\put(0,10){\line(1,-1){10}}
\put(0,10){\line(1,0){10}}
\put(10,10){\line(-1,-1){10}}
\put(0,0){\line(0,1){10}}
\put(3.5,-3){$G_3$}
\end{picture}
\begin{picture}(20,12)(0,0)
\put(5,5){\circle*{2}}
\put(0,0){\circle*{2}}
\put(10,0){\circle*{2}}
\put(10,10){\circle*{2}}
\put(0,10){\circle*{2}}
\put(0,10){\line(1,-1){10}}
\put(0,10){\line(1,0){10}}
\put(10,10){\line(-1,-1){10}}
\put(10,0){\line(0,1){10}}
\put(0,0){\line(0,1){10}}
\put(3.5,-3){$G_4$}
\end{picture}
\begin{picture}(20,12)(0,0)
\put(5,5){\circle*{2}}
\put(0,0){\circle*{2}}
\put(10,0){\circle*{2}}
\put(10,10){\circle*{2}}
\put(0,10){\circle*{2}}
\put(0,10){\line(1,-1){10}}
\put(0,10){\line(1,0){10}}
\put(10,10){\line(-1,-1){10}}
\put(10,0){\line(0,1){10}}
\put(0,0){\line(0,1){10}}
\put(0,0){\line(1,0){10}}
\put(3.5,-3){$G_5$}
\end{picture}
\]
\[
\ \ \ \ \ \
\begin{picture}(23,12)(3,5)
\put(15,0){\circle*{2}}
\put(0,0){\circle*{2}}
\put(7.5,0){\circle*{2}}
\put(0,10){\circle*{2}}
\put(7.5,10){\circle*{2}}
\put(15,10){\circle*{2}}
\put(15,0){\line(0,1){10}}
%
\put(15,10){\line(-3,-4){7.5}}
\put(7.5,0){\line(-3,4){7.5}}
\put(7.5,0){\line(0,1){10}}
\put(0,0){\line(0,1){10}}
\put(0,0){\line(1,0){15}}
%
\put(3,-3){$G_6$}
\end{picture}
\ \ \
\begin{picture}(23,12)(2,5)
\put(15,0){\circle*{2}}
\put(0,0){\circle*{2}}
\put(7.5,0){\circle*{2}}
\put(0,10){\circle*{2}}
\put(7.5,10){\circle*{2}}
\put(15,10){\circle*{2}}
\put(0,0){\line(3,4){7.5}}
\put(15,0){\line(0,1){10}}
\put(15,0){\line(-3,2){15}}
\put(15,0){\line(-3,4){7.5}}
\put(15,10){\line(-3,-2){15}}
\put(15,10){\line(-3,-4){7.5}}
\put(7.5,0){\line(-3,4){7.5}}
\put(7.5,0){\line(0,1){10}}
\put(0,0){\line(0,1){10}}
\put(0,0){\line(1,0){7.5}}
\put(0,10){\line(1,0){7.5}}
\put(3,-3){$G_7$}
\end{picture}
\ \ \
\begin{picture}(22,12)(1,0)
\put(5,0){\circle*{2}}
\put(10,0){\circle*{2}}
\put(15,0){\circle*{2}}
\put(10,-5){\circle*{2}}
\put(10,5){\circle*{2}}
\put(0,0){\circle*{2}}
\put(0,0){\line(1,0){15}}
\put(5,0){\line(1,1){5}}
\put(5,0){\line(1,-1){5}}
\put(15,0){\line(-1,1){5}}
\put(15,0){\line(-1,-1){5}}
\put(10,0){\line(0,1){5}}
\put(2,-8){$G_{8}$}
\end{picture}
\ \ \
\begin{picture}(22,12)(0,0)
\put(5,0){\circle*{2}}
\put(10,0){\circle*{2}}
\put(15,0){\circle*{2}}
\put(15,-5){\circle*{2}}
\put(15,5){\circle*{2}}
\put(0,0){\circle*{2}}
\put(0,0){\line(1,0){15}}
\put(10,0){\line(1,1){5}}
\put(10,0){\line(1,-1){5}}
\put(15,0){\line(0,1){5}}
\put(6,-8){$G_{9}$}
\end{picture}
\]
\[
\ \ \ \ \
\begin{picture}(29,20)(0,-2)
\put(5,0){\circle*{2}}
\put(10,0){\circle*{2}}
\put(15,0){\circle*{2}}
\put(20,-5){\circle*{2}}
\put(20,5){\circle*{2}}
\put(0,0){\circle*{2}}
\put(0,0){\line(1,0){15}}
\put(15,0){\line(1,1){5}}
\put(15,0){\line(1,-1){5}}
\put(7.5,-5.5){$G_{10}$}
\end{picture}
\ \ \
\begin{picture}(29,20)(0,-2)
\put(5,0){\circle*{2}}
\put(10,0){\circle*{2}}
\put(15,5){\circle*{2}}
\put(15,-5){\circle*{2}}
\put(20,0){\circle*{2}}
\put(0,0){\circle*{2}}
\put(0,0){\line(1,0){10}}
\put(10,0){\line(1,1){5}}
\put(20,0){\line(-1,-1){5}}
\put(20,0){\line(-1,1){5}}
\put(10,0){\line(1,-1){5}}
\put(6.5,-5.5){$G_{11}$}
\end{picture}
\ \ \
\begin{picture}(29,20)(0,-2)
\put(5,0){\circle*{2}}
\put(10,0){\circle*{2}}
\put(15,0){\circle*{2}}
\put(20,-5){\circle*{2}}
\put(20,5){\circle*{2}}
\put(0,0){\circle*{2}}
\put(0,0){\line(1,0){15}}
\put(15,0){\line(1,1){5}}
\put(15,0){\line(1,-1){5}}
\put(20,-5){\line(0,1){10}}
\put(7,-5.5){$G_{12}$}
\end{picture}
\ \ \
\begin{picture}(29,20)(0,-2)
\put(5,0){\circle*{2}}
\put(10,5){\circle*{2}}
\put(15,0){\circle*{2}}
\put(10,-5){\circle*{2}}
\put(20,0){\circle*{2}}
\put(0,0){\circle*{2}}
\put(0,0){\line(1,0){5}}
\put(5,0){\line(1,1){5}}
\put(5,0){\line(1,-1){5}}
\put(15,0){\line(1,0){5}}
\put(15,0){\line(-1,1){5}}
\put(15,0){\line(-1,-1){5}}
\put(10,-5){\line(0,1){10}}
\put(2,-5.5){$G_{13}$}
\end{picture}
\]
\end{lem}
~\\[5pt]
{\bf Proof.}
Graphs $G_1,\ldots,G_5$, $G_{10}$ and $G_{11}$ have two eigenvalues less than $-1$,
and the other ones have three positive eigenvalues.
Therefore, by interlacing, none of these graphs is an induced subgraph of a graph in $\ggg$.
\qed
\section{Bipartite graphs}\label{bip}
It is well known that the spectrum of a bipartite graph is symmetric around $0$.
Therefore the bipartite graphs in $\ggg$ are precisely the connected bipartite graphs with
just two positive eigenvalues, one of which equals $1$.
The next result appeared in \cite{fq} (see also \cite{o}).
\begin{lem}\label{bip1}
If $H$ is a connected bipartite graph with exactly two positive eigenvalues,
then $H$ is a mixed coclique extension of the path $P_4$ or $P_5$.
\end{lem}
{\bf Proof.}
Let $B$ be the adjacency matrix of $H$.
Then we can take
\[
B=\left[\begin{array}{cc}
O & N \\
N^\top & O
\end{array}\right].
\]
Since $B$ has exactly two positive eigenvalues,
$B$ also has exactly two negative eigenvalues, therefore $B$ has rank~$4$, and $N$ has rank~$2$.
This leads to just two possible structures for $N$ or $N^\top$:
\[
\left[\begin{array}{cc}
J & O \\
J & J
\end{array}\right],\ {\rm or}
\ \ \left[\begin{array}{ccc}
J & J & O\\
J & O & J
\end{array}\right].
\]
In the first case $H$ is a mixed coclique extension of $P_4$,
in the second case $H$ is a mixed coclique extension of $P_5$.
\qed

\begin{thm}\label{bip2}
If $H$ is a bipartite graph in $\ggg$, then $H$ is one of the following.
\\
(i) A mixed extensiom of $P_4$ of type $(-2,\ ,\ ,-2)$, $(\ ,-2,\ , -3)$, or $(\ ,-3,-2,-2)$,
\\
(ii) A mixed extension of $P_5$ of type $(\ ,\ ,-r,\ ,\ )$ with $r\geq 1$.
\end{thm}
{\bf Proof.}
We have to determine which graphs from Lemma~\ref{bip1} have an eigenvalue $1$,
therefore we need to check when the quotient matrix $Q$ has an eigenvalue $1$.
If $H$ is a mixed extension of $P_4$ of type $(-p,-q,-r,-s)$, then
\[
Q=\left[
\begin{array}{cc}
O & R\\
S & O
\end{array}
\right],
\ {\rm with}
\ \ R=\left[
\begin{array}{cc}
q & 0\\
q & s
\end{array}
\right],\ {\rm and}
\ \ S=\left[
\begin{array}{cc}
p & r\\
0 & r
\end{array}
\right].
\]
The matrix $Q$ has an eigenvalue $1$ whenever $Q^2$ has an eigenvalue $1$, which is the case if and only if $RS$ has an eigenvalue $1$.
We easily have that $\det(RS-I)=pqrs-pq-qr-rs+1$.
So we need to solve $pqrs-pq-qr-rs+1=0$.
We rewrite this equation as $pq(rs-3) + qr(ps-3) + rs(pq-3) = -3$.
Since at least one term on the left hand side is negative we have that $rs\leq 2$, $ps\leq 2$, or $pq\leq 2$.
Each of these three cases give three possibilities, each of which is easily worked out.
Taking the symmetry of the path into account, this leads to the three solutions given in (i).

Next suppose $H$ is a mixed extension of $P_5$.
Because of forbidden induced subgraphs $G_{10}$ and $G_{11}$, the mixed extension can only be of type $(\ ,\ ,-r,\ ,\ )$,
and for all $r\geq 1$ the corresponding quotient matrix $Q$ has an eigenvalue $-1$, and an eigenvalue $0$.
\qed
\section{Reduction}
We call a graph $G$ {\em reduced} if $G$ is not a mixed extension of a smaller graph.
It is tempting to believe that every graph is a mixed extension of a reduced graph.
However, this is not true.
For example the path $P_3$ is a mixed extension of $K_2$ (and not of $K_1$),
but $K_2$ is not reduced.
Because of this we must distinguish mixed clique extensions and mixed coclique extensions.
The next theorem is the key result for the determination of $\ggg$.
\begin{thm}\label{key}
Every graph in $\ggg$ is a mixed clique extension of a connected bipartite graph with at most two positive eigenvalues.
\end{thm}
{\bf Proof.}
We call an edge $\{x,y\}$ of a graph $G$ {\em reducible} if $G$ contains no vertex, other than $x$ and $y$,
that is adjacent to exactly one vertex of $\{x,y\}$.
In terms of the adjacency matrix $A$, being reducible means that the two rows of $A+I$ corresponding $x$ and $y$ are identical.
Let $H\in\ggg$.
We can assume that $H$ is a mixed clique extension of a graph $G$ with no reducible edges.
Clearly $G$ is an induced subgraph of $H$, and therefore $G$ contains none of the forbidden graphs
given in Section~\ref{forbidden}.
In particular, $G$ contains no odd cycle $C_n$ with $n\geq 5$.

Assume $G$ contains a triangle $\T$.
If $G$ has a vertex not in $\T$, but adjacent to exactly two vertices $x$ and $y$ (say) of $\T$.
Then any other vertex not in $\T$ is either adjacent to $x$ and $y$ or not adjacent to $x$ and $y$,
since otherwise $G$ has $G_3$, $G_4$, or $G_5$ as an induced subgraph.
This implies that the edge $\{x,y\}$ is reducible, which is a contradiction.
Similarly, assume $G$ has a vertex nonadjacent to exactly two vertices $x$ and $y$ of $\T$.
Then again it follows that every other vertex is either adjacent to $x$ and $y$ or nonadjacent to $x$ and $y$,
since otherwise $G$ has $G_1$, $G_2$, $G_3$, or $G_4$ as an induced subgraph.
So again we get a contradiction.
It is also not possible that every vertex not in $\T$ is adjacent to all or no vertices of $\T$, since then each edge of $\T$ is reducible.
Therefore $G$ has no triangle, and hence $G$ is bipartite.
\qed
\\

A connected bipartite graph with just one positive eigenvalue is a complete bipartite graph,
therefore each graph $H\in\ggg$ is a mixed clique extension of a complete bipartite graph, or a
mixed clique extension of a graph described in Lemma~\ref{bip1}.

\section{The determination of $\ggg$}

\begin{thm}\label{main}
A graph $H$ belongs to $\ggg$ if and only if $H$ is one of the following.
\\[3pt]
(i)
A mixed extension of $P_3$ of type $(-p,-q,r),\ (-p,q,r),\ (p,-q,r)$, or $(p,q,r)$
\\
\hspace*{12pt} with $p,q\geq 1$ and $r\geq 2$,
\\[3pt]
(ii)
a mixed extension of $P_4$ of type $(p,-3,-2,-2)$, $(-2,q,r,-2)$, or $(p,-2,r,-3)$
\\
\hspace*{12pt} with $p,q,r\geq 1$,
\\[3pt]
(iii)
a mixed extension of $P_4$ of type $(p,q,-r,s)$, with $r\geq 1$ and $(p,q,s)\in$
\[
\{(3,3,6),(3,4,4),(3,6,3),(4,2,6),(4,3,3),(4,6,2),(5,2,4),(5,4,2),(7,2,3),(7,3,2)\},
\]
(iv)
a mixed extension of $P_4$ of type $(p,q,r,s)$, with $(p,q,r,s)\in$
\[
\{(2,2,2,7),(2,2,3,4),(2,2,6,3),(2,3,2,5),(2,3,4,3),(2,5,2,4),(2,5,3,3),(3,2,2,3)\},
\]
(v)
a mixed extension of $P_5$ of type $(\ ,p,-q,r,\ )$ with $p,q,r\geq 1$.
\end{thm}
{\bf Proof.}
By Theorem~\ref{key} and Lemma~\ref{bip1}, $H$ is a mixed clique extension of the complete bipartite graph $K_{p,q}$,
or a mixed clique extension of a mixed coclique extension of $P_4$ or $P_5$.
First assume that $H$ is a mixed clique extension of $K_{p,q}$.
We know $H\neq K_{p,q}$, and the forbidden induced subgraphs $G_6$ (which is a mixed extension of $K_{1,3}$ of type $(\ ,2,2,\ )$) and $G_7$
(which is a mixed extension of $K_{2,2}$ of type $(\ ,2,2,\ )$) lead to the conclusion that $H$ is a mixed extension of $K_{1,2}=P_3$.
Using the symmetry of $P_3$ and excluding the mixed extensions of $K_2$ we obtain the mentioned cases.

Next assume $H$ is a mixed clique extension of a graph $G'$ which in turn is a mixed coclique extension of $P_4$ or $P_5$.
Suppose a vertex $x$ of $P_4$ or $P_5$ is replaced by a coclique $V_x$ of $G'$ in the first step,
and that a vertex $x'\in V_x$ is replaced by a clique $V'_{x'}$, in the second step.
Then $|V_x|=1$ or $|V'_{x'}|=1$, because otherwise $H$ contains the forbidden induced subgraph $G_8$ or $G_9$.
This implies that $H$ is a mixed extension of $P_4$ or $P_5$.

First we consider $P_4$, and examine all possible types: $(\pm p,\pm q,\pm r,\pm s)$ with $p,q,r,s \geq 1$.
We investigate each of the sixteen possible sign patterns, but because of symmetry, and because the bipartite case (all signs negative)
is already done in Theorem~\ref{bip2}, there are nine cases to be checked.
For each case we consider the quotient matrix $Q$ and check if
$Q$ has an eigenvalue $-1$ or $0$, which is the case if and only if $H\in\ggg$.
\\[3pt]
Type $(p,q,r,-s)$:  $\det(Q)=rs(p+q-1)\neq 0$, $\det(Q+I)=-pqr\neq 0$. No solution.
\\[3pt]
Type $(p,q,-r,-s)$: $\det(Q)=rs(p+q-1)\neq 0$, $\det(Q+I)=-pqr\neq 0$. No solution.
\\[3pt]
Type $(p,-q,-r,s)$: $\det(Q)=qr(p+s-1)\neq 0$, $\det(Q+I)=ps(1-q-r)\neq 0$. No solution.
\\[3pt]
Type $(p,-q,r,-s)$: $\det(Q)=pqrs\neq 0$, $\det(Q+I)=pr(qs-s-2q+1)$.
It follows that $Q$ has eigenvalue $-1$ if and only if $s(q-2)+q(s-4)=-2$.
Since one of the two terms of the left hand side is negative, it follows that $q<2$, or $s<4$.
Using this, it follows straightforwardly that $s=3$, $q=2$ is the only solution.
Thus we find type $(p,-2,r,-3)$ of item (ii).
\\[3pt]
Type $(p,-q,-r,-s)$: $\det(Q)=pqrs\neq 0$, $\det(Q+I)=p(qrs-rs-qr-q+1)$.
So $Q$ has an eigenvalue $-1$ if $rs(q-3)+q(rs-3)+qr(s-3)=-3$.
This gives that $q\leq 2$, $rs\leq 2$ or $s\leq 2$.
If $r=1$ we deal with the special case $(p,-2,\ ,-3)$ of the previous case.
If $r\geq 2$ we find $q=3$, $r=s=2$, which gives type $(p,-3,-2,-2)$ of item (ii).
\\[3pt]
Type $(-p,q,r,-s)$: $\det(Q)=pqrs\neq 0$, $\det(Q+I)=qr((p-1)(s-1)-1)$.
So $\det(Q+I)=0$ if $p=s=2$.
This gives type $(-2,q,r,-2)$ of item (ii).
\\[3pt]
Type $(-p,-q,r,-s)$: $\det(Q)=pqrs\neq 0$, $\det(Q+I)=r(pqs-pq-q-s+1)$.
If $p=1$, $q=1$, or $s=1$ we deal with a special case of one of the previous types
($(\ ,-2,r,-3)$ is the solution if $p=1$ and $(-2,\ ,r,-2)$ is the solution if $q=1$).
So we can assume that $p,q,s\geq 2$, which implies $pqs-pq-q-s\geq 0$, so there is no
new solution for this type.
\\[3pt]
Type $(p,q,-r,s)$: $\det(Q+I)=-pqrs\neq 0$, $\det(Q)=r(-pqs+2qs +pq+ps-s-q)$.
So $Q$ has an eigenvalue $0$ if and only if $qs(p-6) + pq(s-3) + ps(q-3)=-3(q+s)$.
This implies that $p\leq 5$, $s\leq 2$, or $q\leq 2$.
The cases $p=1$, $q=1$ and $s=1$ have been considered already, so we may assume $p,q,s\geq 2$.
Working out each of the possibilities for $p$, $q$ and $s$ we find the solutions given in item (iii).
\\[3pt]
Type $(p,q,r,s)$: $\det(Q+I)=-pqrs\neq 0$, $\det(Q)=-pqrs+pqr+qrs+ps+pr+qs-p-q-r-s+1$.
So $\det(Q)=0$ if and only if
$3pqr(s-4) + 3qrs(p-4) + 2ps(qr-6) + 2pr(qs-6) + 2qs(pr-6)=-12(p+q+r+s-1)$.
This gives $s\leq 3$, $p\leq 3$, $qr\leq 5$, $qs\leq 5$, or $pr\leq 5$.
Also because all other types have been considered before, we may assume that $p,q,r,s\geq 2$.
Because of the symmetry of $P_4$ we may also assume that $p\leq s$, and $q\leq r$ if $p=s$.
Working out all possibilities (which is in this case a tedious job) leads to the solutions of item (iv).

Finally suppose $H$ is a mixed extension of $P_5$.
The forbidden subgraphs $G_{10},\ldots,G_{13}$ imply that the type can only be $(\ ,p,-q,r,\ )$.
But then the quotient matrix $Q$ has an eigenvalue $-1$ and an eigenvalue $0$ for all $p,q,r\geq 1$.
This proves (v).
\qed\\

Note that Theorem~\ref{bip2} is a special case of Theorem~\ref{main}.
Thus the complete determination of $G'$ is given by Theorem~\ref{main} and the three propositions of Section~\ref{easy}.
Combining these results we obtain:
\begin{cor}
A connected graph $H$ has all but at most three eigenvalues equal to $-1$ or $0$ if and only if $H$ is a mixed extension
of a connected graph with no more than three vertices, or $H$ is a mixed extension of $P_4$ or $P_5$ of one of the types given in Theorem~\ref{main}(ii)-(v).
\end{cor}


\end{document}